\documentclass[11pt]{amsart}
\usepackage{amssymb, adjustbox, enumerate, amsbsy, stmaryrd}
\usepackage{array, longtable}
\usepackage{geometry}
\usepackage{todonotes}

\usepackage{amsfonts, amssymb, amscd}
\numberwithin{equation}{section}

\usepackage[symbol]{footmisc}

\usepackage{bm}
\usepackage{verbatim}
\usepackage{mathrsfs}
\usepackage{graphicx}
\usepackage{tikz-cd}
\usepackage{subcaption}
\usepackage{listings}
\usepackage{subfiles}
\usepackage[toc,page]{appendix}
\usepackage{mathtools}
\usepackage{comment}
\usepackage{enumerate}
\usepackage{enumitem}
\usepackage[all]{xy}

\usepackage{graphicx}
\graphicspath{{images/}}

\usepackage{appendix}
\usepackage{hyperref}
\hypersetup{
    colorlinks=true,
    citecolor=red,
    linkcolor=blue,
    filecolor=magenta,      
    urlcolor=red,
}
\newcommand{\bb}{\bm{b}}
\newcommand{\Mm}{{\bf{M}}}

\newcommand{\Dd}{{\bf{D}}}

\newcommand{\Spec}{\mathrm{Spec}}

\newcommand{\Qq}{\mathbb{Q}}

\newcommand{\Rr}{\mathbb{R}}

\newcommand{\Center}{\operatorname{center}}

\newcommand{\Hom}{\operatorname{Hom}}

\newcommand{\Ker}{\operatorname{Ker}}
\newcommand{\Ima}{\operatorname{Im}}

\newcommand{\Nklt}{\operatorname{Nklt}}

\newcommand{\Supp}{\operatorname{Supp}}
\newcommand{\Bs}{\operatorname{Bs}}

\newcommand{\Nlc}{\operatorname{Nlc}}

\newcommand{\mult}{\operatorname{mult}}

\newcommand{\cont}{\operatorname{cont}}

\newcommand{\Oo}{\mathcal{O}}

\newtheorem{thm}{Theorem}[section]

\newtheorem{lem}[thm]{Lemma}

\newtheorem{claim}[thm]{Claim}

\theoremstyle{definition}
\newtheorem{defn}[thm]{Definition}

\theoremstyle{definition}
\newtheorem{rem}[thm]{Remark}

\theoremstyle{definition}

\begin{document}

\title{Vanishing theorems for generalized pairs}

\author{Bingyi Chen, Jihao Liu, and Lingyao Xie}

\address{Yau Mathematical Sciences Center, Tsinghua University, Beijing 100084, P. R. China}
\email{bychen@mail.tsinghua.edu.cn}

\address{Department of Mathematics, Northwestern University, 2033 Sheridan Rd, Evanston, IL 60208, USA}
\email{jliu@northwestern.edu}

\address{Department of Mathematics, The University of Utah, 155 South 1400 East, JWB 233, Salt Lake City, UT 84112, USA}
\email{lingyao@math.utah.edu}

\subjclass[2020]{14E30, 14B05}
\keywords{Vanishing theorem. Generalized pairs. Lc singularity.}
\date{\today}

\begin{abstract}
We establish the Kodaira vanishing theorem and the Kawamata-Viehweg vanishing theorem for lc generalized pairs. As a consequence, we provide a new proof of the base-point-freeness theorem for lc generalized pairs. This new approach allows us to prove the contraction theorem for lc generalized pairs without using Koll\'ar's gluing theory.
\end{abstract}

\maketitle
\tableofcontents

\section{Introduction}

We work with the field of complex numbers $\mathbb{C}$. All generalized pairs are assumed to be NQC generalized pairs (cf. \cite{HL22}) in this paper.

The theory of ``generalized pairs" (abbreviated as ``g-pairs") holds significant importance in modern birational geometry. It was initially introduced by Birkar and Zhang in their study on effective Iitaka fibrations \cite{BZ16}. Since then, this theory has proven to be crucial in various aspects of birational geometry, including the proof of the Borisov-Alexeev-Borisov conjecture \cite{Bir19, Bir21a}, the theory of complements \cite{Bir19, Sho20}, the connectedness principle \cite{Bir20, FS23}, non-vanishing theorems \cite{LMPTX22}, the minimal model program for K\"ahler manifolds \cite{DHY23, HX23}, and foliations \cite{LLM23}, etc. For a comprehensive overview of the theory of g-pairs, we refer interested readers to \cite{Bir21b}.

An important aspect of the study of g-pairs is their minimal model program. The foundations for the minimal model program of klt g-pairs and $\mathbb{Q}$-factorial dlt g-pairs were established in \cite{BZ16,HL22}. Recently, progress has been made towards the minimal model program theory for lc g-pairs. Specifically, a series of recent works \cite{HL21, LX22, Xie22} have established the cone theorem, contraction theorem, base-point-freeness theorem, and the existence of flips for lc g-pairs. This enables us to run the minimal model program for lc g-pairs in a comprehensive manner. For further related works, we refer readers to \cite{Has22a,Has22b,LT22,LX23,TX23}.

Apart from the minimal model program, there are numerous other topics within classical birational geometry that are worth discussing in the context of lc g-pairs. For instance, it is known that lc g-pairs have Du Bois singularities \cite{LX22}. In this paper, we establish several vanishing theorems for lc g-pairs. The first main theorem of the paper is the following:

\begin{thm}\label{thm: kod vanishing with lc strata}
Let $(X,B,\Mm)/U$ be an lc generalized pair associated with projective morphism $f: X\rightarrow U$, $D$ a Cartier divisor on $X$ such that $D-(K_X+B+\Mm_X)$ is nef$/U$ and log big$/U$ with respect to $(X,B,\Mm)$ (cf. Definition \ref{defn: log big}), $Y$ a union of lc centers of $(X,B,\Mm)$ such that $Y\not=X$, and $\mathcal{I}_Y$ the defining ideal sheaf of $Y$ on $X$. Then:
\begin{enumerate}
	\item $R^if_*\mathcal{O}_Y(D)=0$ for any positive integer $i$.
	\item $R^if_*\mathcal{O}_X(D)=0$ for any positive integer $i$.
    \item The map $f_*\mathcal{O}_X(D)\rightarrow f_*\mathcal{O}_Y(D)$ is surjective.
    \item $R^if_*(\mathcal{I}_Y\otimes\mathcal{O}_X(D))=0$ for any positive integer $i$.
\end{enumerate}
\end{thm}

\begin{rem}
We briefly explain the history on results that are related to Theorem \ref{thm: kod vanishing with lc strata}.
\begin{enumerate}
    \item When $\Mm=\bm{0}$ and $(X,B)$ is klt, Theorem \ref{thm: kod vanishing with lc strata}(1)(3) become trivial, and Theorem \ref{thm: kod vanishing with lc strata}(2)(4) are both equivalent to the usual relative Kawamata-Viehweg vanishing theorem (cf. \cite[Theorem 1-2-7]{KMM87}).
    \item When $\Mm=\bm{0}$ and $D-(K_X+B+\Mm_X)$ is ample$/U$, Theorem \ref{thm: kod vanishing with lc strata}(2) becomes the usual Kodaira vanishing theorem for lc pairs \cite[Theorem 4.4]{Fuj09} and Theorem \ref{thm: kod vanishing with lc strata}(4) is \cite[Theorem 7.3]{Amb03} and  \cite[Theorem 8.1]{Fuj11}. 
    \item When $\Mm=\bm{0}$, Theorem \ref{thm: kod vanishing with lc strata}(2)(3)(4) follow from \cite[Theorem 7.3]{Amb03}, \cite[Theorem 6.3.4(2)]{Fuj17} and  Theorem \ref{thm: kod vanishing with lc strata}(1) follows from \cite[Theorem 1.14]{Fuj14}. Note that Theorem \ref{thm: kod vanishing with lc strata}(2) becomes the usual Kawamata-Viehweg vanishing theorem for lc pairs. 
    \item In fact, \cite[Theorem 7.3]{Amb03}, \cite[Theorem 6.3.4(2)]{Fuj17} prove the qlc case of Theorem \ref{thm: kod vanishing with lc strata}. Since any qlc pair is always an lc g-pair \cite[Remark 1.9]{Fuj22}, Theorem \ref{thm: kod vanishing with lc strata} implies  \cite[Theorem 7.3]{Amb03}, \cite[Theorem 6.3.4(2)]{Fuj17} for qlc pairs.
    \item There is no previously written result when $\Mm\not=\bm{0}$, but the case when $\Mm_X$ is $\Rr$-Cartier and $D-(K_X+B+\Mm_X)$ is ample$/U$ can be easily deduced from \cite[Lemma 5.18]{HL21} and the Kodaira vanishing theorem for lc pairs. 
\end{enumerate}
\end{rem}

Theorem \ref{thm: kod vanishing with lc strata} immediately implies the Kodaira vanishing theorem for lc g-pairs and the Kawamata-Viehweg vanishing theorem for lc g-pairs. We provide the precise statement of these results here as they are more useful for direct applications.

\begin{thm}[Kodaira vanishing theorem for lc generalized pairs]\label{thm: kod vanishing gpair intro}
Let $(X,B,\Mm)$ be a projective lc generalized pair, and let $D$ be a Cartier divisor on $X$ such that $D-(K_X+B+\Mm_X)$ is ample. Then $H^i(X,\mathcal{O}_X(D))=0$ for any positive integer $i$.
\end{thm}

\begin{thm}[Relative Kawamata-Viehweg vanishing for lc generalized pairs]\label{thm: kv vanishing gpair intro}
Let $(X,B,\Mm)/U$ be an lc generalized pair associated with morphism $f: X\rightarrow U$, and let $D$ be a Cartier divisor on $X$ such that $D-(K_X+B+\Mm_X)$ is nef$/U$ and log big$/U$ with respect to $(X,B,\Mm)$. Then $R^if_*\mathcal{O}_X(D)=0$ for any positive integer $i$.
\end{thm}

It was anticipated by Hashizume \cite[Page 77, Line 24-25]{Has22b} that Theorem \ref{thm: kod vanishing gpair intro} would play a pivital role in establishing the base-point-freeness theorem for lc g-pairs. Despite the base-point-freeness theorem's prior proof in \cite[Theorem 1.4]{Xie22}, we endeavor to explore the viability of Hashizume's approach. Leveraging the implications of Theorem \ref{thm: kod vanishing with lc strata}, we provide a new proof of the base-point-free theorem for lc g-pairs, thereby fulfilling Hashizume's expectation. It is noteworthy that our proof diverges significantly from the one in \cite{Xie22}, as the latter relies heavily on Koll\'ar's gluing theory for g-pairs, while our novel approach bypasses this necessity.

\begin{thm}[Base-point-freeness theorem for lc generalized pairs, cf. {\cite[Theorem 1.4]{Xie22}}]\label{thm:base-point-freeness intro}
Let $(X,B,\Mm)/U$ be an lc g-pair and $D$ a nef$/U$ Cartier divisor on $X$, such that $aD-(K_X+B+\Mm_X)$ is ample$/U$ for some positive real number $a$. Then $\mathcal{O}_X(mD)$ is globally generated over $U$ for any integer $m\gg 0$.
\end{thm}

As an immediate application, we have the following semi-ampleness theorem for lc g-pairs.

\begin{thm}[Semi-ampleness theorem for lc generalized pairs, cf. {\cite[Theorems 1.2]{Xie22}}]\label{thm: semi-ampleness intro}
Let $(X,B,\Mm)/U$ be an lc g-pair and $D$ a nef$/U$ $\mathbb R$-Cartier $\Rr$-divisor on $X$, such that $D-(K_X+B+\Mm_X)$ is ample$/U$. Then $D$ is semi-ample$/U$.
\end{thm}

We remark that \cite[Theorem 1.4]{Xie22} is stronger than Theorem \ref{thm:base-point-freeness intro} since \cite[Theorem 1.4]{Xie22} only requires that $aD-(K_X+B+\Mm_X)$ is nef$/U$ and log big$/U$. Nonetheless, Theorem \ref{thm:base-point-freeness intro} is strong enough for us to immediately deduce the contraction theorem for lc g-pairs \cite[Theorem 1.5]{Xie22} without using Koll\'ar's gluing theory (see Remark \ref{rem: does not rely on kollar}).

\begin{thm}[Contraction theorem for lc generalized pairs, cf. {\cite[Theorem 1.5]{Xie22}}]\label{thm: cont thm gpair}
Let $(X,B,\Mm)/U$ be an lc generalized pair and $F$ a $(K_X+B+\Mm_X)$-negative extremal face$/U$. Then there exists a contraction$/U$ $\cont_F: X\rightarrow Z$ of $F$ satisfying the following.
\begin{enumerate}
    \item For any integral curve $C$ on $X$ such that the image of $C$ in $U$ is a closed point, $\cont_F(C)$ is a point if and only if $[C]\in F$.
    \item $\mathcal{O}_Y=(\cont_F)_*\mathcal{O}_X$. In other words, $\cont_F$ is a contraction.
    \item For any Cartier divisor $D$ on $Y$ such that $D\cdot C=0$ for any curve $C$ contracted by $\cont_F$, there exists a Cartier divisor $D_Y$ on $Y$ such that $D=\cont_F^*D_Y$.
\end{enumerate}
\end{thm}

\medskip

\noindent\textbf{Acknowledgement}. The authors would like to thank Christopher D. Hacon for valuable discussions and constant support to the second and third authors. The first author would like to thank Caucher Birkar for constant support. The second author would like to thank Yuchen Liu for constant support. The project started when Caucher Birkar kindly invited the second author to visit Tsinghua University in April 2023. The second author would like to extend their appreciation for the warm hospitality received during the visit. A portion of the project was carried out during the first and second authors' visit to Fudan University in April 2023, and the authors would like to express their thanks for the hospitality extended to them. The authors would like to thank Jingjun Han and Fanjun Meng for useful discussions. We would like to acknowledge the assistance of ChatGPT in polishing the wording. The third author is partially supported by NSF research grants no: DMS-1801851, DMS-1952522 and by a grant from the Simons Foundation; Award Number: 256202.

\section{Preliminaries}

Throughout the paper, we will mainly work with normal quasi-projective varieties to ensure consistency with the references. However, most results should also hold for normal varieties that are not necessarily quasi-projective. Similarly, most results in our paper should hold for any algebraically closed field of characteristic zero. We will adopt the standard notations and definitions in \cite{KM98, BCHM10} and use them freely. For generalized pairs, we will follow the notations as in \cite{HL21}.

\subsection{Definition of generalized pairs}

\begin{defn}[$\bb$-divisors]\label{defn: b divisors} Let $X$ be a normal quasi-projective variety. We call $Y$ a \emph{birational model} over $X$ if there exists a projective birational morphism $Y\to X$. 

Let $X\dashrightarrow X'$ be a birational map. For any valuation $\nu$ over $X$, we define $\nu_{X'}$ to be the center of $\nu$ on $X'$. A \emph{$\bb$-divisor} $\Dd$ over $X$ is a formal sum $\Dd=\sum_{\nu} r_{\nu}\nu$ where $\nu$ are valuations over $X$ and $r_{\nu}\in\mathbb R$, such that $\nu_X$ is not a divisor except for finitely many $\nu$. The \emph{trace} of $\Dd$ on $X'$ is the $\Rr$-divisor
$$\Dd_{X'}:=\sum_{\nu_{X'}\text{ is a divisor}}r_\nu\nu_{X'}.$$
If $\Dd_{X'}$ is $\Rr$-Cartier and $\Dd_{Y}$ is the pullback of $\Dd_{X'}$ on $Y$ for any birational model $Y$ over $X'$, we say that $\Dd$ \emph{descends} to $X'$ and $\Dd$ is the \emph{closure} of $\Dd_{X'}$, and write $\Dd=\overline{\Dd_{X'}}$. 

Let $X\rightarrow U$ be a projective morphism and assume that $\Dd$ is a $\bb$-divisor over $X$ such that $\Dd$ descends to a birational model $Y$ over $X$. If $\Dd_Y$ is nef$/U$, then we say that $\Dd$ is \emph{nef}$/U$. If $\Dd_Y$ is a Cartier divisor, then we say that $\Dd$ is \emph{$\bb$-Cartier}. If $\Dd$ can be written as an $\Rr_{\geq 0}$-linear combination of nef$/U$ $\bb$-Cartier $\bb$-divisors, then we say that $\Dd$ is \emph{NQC}$/U$.

We let $\bm{0}$ be the $\bb$-divisor $\bar{0}$.
\end{defn}

\begin{defn}[Generalized pairs]\label{defn: g-pairs}
A \emph{generalized pair} (\emph{g-pair} for short) $(X,B,\Mm)/U$ consists of a normal quasi-projective variety $X$ associated with a projective morphism $X\rightarrow U$, an $\Rr$-divisor $B\geq 0$ on $X$, and an NQC$/U$ $\bb$-divisor $\Mm$ over $X$, such that $K_X+B+\Mm_X$ is $\Rr$-Cartier.

If $\Mm=\bf{0}$, a g-pair $(X,B,\Mm)/U$ is called a \emph{pair} and is denoted by $(X,B)$ or $(X,B)/U$. 

If $U=\{pt\}$, we usually drop $U$ and say that $(X,B,\Mm)$ is \emph{projective}. If $U$ is not important, we may also drop $U$.
\end{defn}

\begin{defn}[Singularities of generalized pairs]\label{defn: sing of g-pairs}
	Let $(X,B,\Mm)/U$ be a g-pair. For any prime divisor $E$ and $\mathbb R$-divisor $D$ on $X$, we define $\mult_{E}D$ to be the \emph{multiplicity} of $E$ along $D$.  Let $h: W\to X$
	be any log resolution of $(X,\Supp B)$ such that $\Mm$ descends to $W$, and let
	$$K_W+B_W+\Mm_W:=h^*(K_X+B+\Mm_X).$$
	The \emph{log discrepancy} of a prime divisor $D$ on $W$ with respect to $(X,B,\Mm)$ is $1-\mult_{D}B_W$ and it is denoted by $a(D,X,B,\Mm).$
	
	We say that $(X,B,\Mm)$ is \emph{lc} (resp. \emph{klt}) if $a(D,X,B,\Mm)\ge0$ (resp. $>0$) for every log resolution $h: W\to X$ as above and every prime divisor $D$ on $W$. We say that $(X,B,\Mm)$ is \emph{dlt} if $(X,B,\Mm)$ is lc, and there exists a closed subset $V\subset X$, such that
\begin{enumerate}
    \item $X\backslash V$ is smooth and $B_{X\backslash V}$ is simple normal crossing, and
    \item for any prime divisor $E$ over $X$ such that $a(E,X,B,\Mm)=0$, $\Center_XE\not\subset V$ and $\Center_XE\backslash V$ is an lc center of $(X\backslash V,B|_{X\backslash V})$.
\end{enumerate}
We refer the reader to \cite[Theorem 6.1]{Has22a} for equivalent definitions of dlt g-pairs. 

Suppose that $(X,B,\Mm)$ is lc. An \emph{lc place} of $(X,B,\Mm)$ is a prime divisor $E$ over $X$ such that $a(E,X,B,\Mm)=0$. An \emph{lc center} of $(X,B,\Mm)$ is either $X$, or the center of an lc place of $(X,B,\Mm)$ on $X$. The \emph{non-klt locus} $\Nklt(X,B,\Mm)$ of $(X,B,\Mm)$ is the union of all lc centers of $(X,B,\Mm)$ except $X$ itself. 

We note that the definitions above are independent of the choice of $U$
\end{defn}

\begin{defn}[Log big]\label{defn: log big}
Let $(X,B,\Mm)/U$ be a g-pair and $D$ an $\Rr$-Cartier $\Rr$-divisor $D$ on $X$. We say that $D$ is \emph{log big$/U$ with respect to $(X,B,\Mm)$} if $D|_V$ is big$/U$ for any lc center $V$ of $(X,B,\Mm)$. In particular, $D$ is big$/U$.
\end{defn}

\subsection{Universal push-out diagram}
\begin{defn}
We say a commutative diagram of schemes
\begin{align*}
\xymatrix{
\mathcal{C}\ar@{->}[r]^j\ar@{->}[d]_q & Y\ar@{->}[d]^{p}\\
 \mathcal{D}\ar@{->}[r]^i & X\\
}
\end{align*}
is a \emph{universal push-out diagram} if for any scheme $T$, the induced diagram
$$\xymatrix{
\Hom{(X,T)}\ar@{->}[r]^{\circ i}\ar@{->}[d]_{\circ p} & \Hom{(\mathcal{D},T)}\ar@{->}[d]^{{\circ q}}\\
 \Hom{(Y,T)}\ar@{->}[r]^{\circ j} & \Hom{(\mathcal{C},T)}\\
}$$
is a universal pull-back diagram of sets. 
\end{defn}
\begin{lem}\label{lem:pushout}
Let $X$ be a semi-normal variety and let $\pi: X^n\to X$ be the normalization of $X$. Let $Z$ be a reduced closed subvariety of $X$ such that $X\setminus Z$ is normal. Let $Y:=\pi^{-1}(Z)$ associated with the reduced scheme structure. Denote the induced morphism $Y\rightarrow Z$ by $\pi_Y$. Then we have the following universal push-out diagram
\begin{center}$\xymatrix{
Y\ar@{->}[d]_{\pi_Y}\ar@{^(->}[r]^{j} & X^n\ar@{->}[d]^{\pi}\\
Z\ar@{^(->}[r]^{i} & X
}$
\end{center}
and a short exact sequence
\begin{align}\label{eq: short exact sequence in upd}
    0\to \Oo_{X}\xrightarrow{\pi^*\oplus i^*} \pi_*\Oo_{X^n}\oplus\Oo_{Z}\xrightarrow{j^*-\pi_Y^*} (\pi_Y)_* \Oo_{Y}\to 0,
\end{align}
where $i,j$ are the natural closed immersions.
\end{lem}
\begin{proof}
Since $j$ is a closed immersion and $\pi_Y$ is a finite morphism, by \cite[Theorem 9.30]{Kol13}, \cite[8.1]{Kol95}, we have a universal push-out diagram
\begin{center}$\xymatrix{
Y\ar@{->}[d]_{\pi_Y}\ar@{^(->}[r]^{j} & X^n\ar@{->}[d]^{\pi'}\\
Z\ar@{^(->}[r]^{i'} & X'
}$
\end{center}
where 
$$
X':=\Spec_X\Ker[\pi_*\Oo_{X^n}\oplus\Oo_{Z}\xrightarrow{j^*-\pi_Y^*} (\pi_Y)_*\Oo_{Y}].
$$
Therefore, it suffices to prove the short exact sequence \eqref{eq: short exact sequence in upd} exists. Let $\mathcal{J}$ be the conductor ideal sheaf of $\pi:X^n\to X$, which can be regarded as both an $\Oo_X$-module and an $\Oo_{X^n}$-module via the inclusion $\Oo_X\hookrightarrow\Oo_{X^n}$. By \cite[5.5.3]{Kol95}, $\mathcal{J}$ is its own radical in $\Oo_{X^n}$ and hence is its own radical in $\Oo_{X}$. Let $\mathcal{I}_{Y},\mathcal{I}_Z$ be the ideal sheaves of $Y,Z$ respectively. Since $X\setminus Z$ is normal, $\mathcal{I}_Z\subset\mathcal{J}$. 
\begin{claim}\label{claim: izox=iz}
$\mathcal{I}_Z\cdot\Oo_{X^n}=\mathcal{I}_Z.$
\end{claim}
\begin{proof}
Let $\mathcal{I}':=\mathcal{I}_Z\cdot\Oo_{X^n}$, then $\mathcal{I}_Z\subset\mathcal{I}'\subset\mathcal{J}$. Since $Z$ is reduced, we only need to prove that the sub-schemes defined by $\mathcal{I}_Z$ and $\mathcal{I}'$ in $X$ have the same support. Thus we only need to prove that the sub-schemes defined by $\mathcal{I}_Z$ and $\mathcal{I}'$ in $X$ have the same support near any point $x\in X$.

If $x\in\Supp \Oo_{X}/\mathcal{J}$, then $x\in\Supp \Oo_{X}/\mathcal{I}'\subseteq Z$ and we are done.

If $x\notin \Supp \Oo_{X}/\mathcal{J}$, then $X$ is normal at $x$, hence $\pi:X^n\to X$ is an isomorphism near $x$. Therefore, $\mathcal{I}_Z=\mathcal{I}'$ and $\Supp \Oo_{X}/\mathcal{I}'=Z$ near $x$. 
\end{proof}

\begin{claim}\label{claim: iz=iy}
$\mathcal{I}_Z=\mathcal{I}_Y$ in $\Oo_{X^n}$. 
\end{claim}
\begin{proof}
By definition, $\mathcal{I}_Y$ is the radical of $\mathcal{I}_Z$ in $\Oo_{X^n}$. Since $\mathcal{I}_Z\subset \mathcal{J}$ and  $\mathcal{J}$ is its own radical in $\Oo_{X^n}$, $\mathcal{I}_Y$ is contained in $\mathcal{J}$ and hence is an ideal sheaf of $\Oo_X$.  Therefore, $\mathcal{I}_Y$ is the radical of $\mathcal{I}_Z$ in $\Oo_X$. Since $Z$ is reduced, $\mathcal{I}_Z=\mathcal{I}_Y$ in $\Oo_{X^n}$.
\end{proof}

\noindent\textit{Proof of Lemma \ref{lem:pushout} continued}. By Claims \ref{claim: izox=iz} and \ref{claim: iz=iy}, we may consider the question locally and assume that $X=\Spec A, X^n=\Spec B$, and $\mathcal{I}_Z=\mathcal{I}_Y=I$. Then the map
$$
\phi: B\oplus A/I\to B/I,~(b,a+I)\mapsto (b-a)+I
$$
is surjective and the map 
$$
\psi: A \to B\oplus A/I,~ a\mapsto (a,a+I)
$$
is injective. Thus
$$(b,a+I)\in\Ker\phi\iff b\in A\text{ and }(b,a+I)=(b,b+I)\iff (b,a+I)\in\Ima(\psi),$$
so \eqref{eq: short exact sequence in upd} is a short exact sequence and we are done.
\end{proof}

\subsection{Union of lc centers of generalized pairs}

\begin{defn}[Union of lc centers] Let $(X,B,\Mm)$ be an lc g-pair. A \emph{union of lc centers} of $(X,B,\Mm)$ is a reduced scheme $Y=\cup Y_i$, where each $Y_i$ is an lc center of $(X,B,\Mm)$. We denote by $S(X,B,\Mm)$ the set of all unions of lc centers of $(X,B,\Mm)$. We remark that 
\begin{enumerate}
    \item $\emptyset$ is also considered as a union of lc centers, and
    \item a union of lc center may be represented in different ways. For example, if $Y_1$ and $Y_2$ are two lc centers such that $Y_1\subsetneq Y_2$, then $Y_1\cup Y_2$ and $Y_2$ are the same union of lc centers.
\end{enumerate}
\end{defn}

\begin{defn}[Adjacent unions of lc centers]
Let $(X,B,\Mm)$ be an lc g-pair. For any two unions of lc centers $Y,Y'\in S(X,B,\Mm)$, we say that $Y$ and $Y'$ are \emph{adjacent} in $S(X,B,\Mm)$ if
\begin{enumerate}
    \item $Y\subsetneq Y'$ or $Y'\subsetneq Y$, and
    \item there does not exist any $Y''\in S(X,B,\Mm)$ such that $Y\subsetneq Y'' \subsetneq Y'$ or $Y'\subsetneq Y''\subsetneq Y$. 
\end{enumerate}
An lc center $V$ is called \emph{minimal} in $S(X,B,\Mm)$ if $V$ and $\emptyset$ are adjacent in $S(X,B,\Mm)$.
\end{defn}

\begin{lem}\label{lem:adjunction}
  Let $(X,B,\Mm)/U$ be an lc g-pair, $W$ a union of lc centers of $(X,B,\Mm)$, and $\pi: W^n\to W$ the normalization of $W$. Suppose that $\dim W\geq 1$. Then there exists an  lc g-pair $(W^n,B_{W^n},\Mm^{W^n})/U$, such that
\begin{enumerate}
    \item $K_{W^n}+B_{W^n}+\Mm^{W^n}_{W^n}\sim_{\Rr,U}(K_X+B+\Mm_X)|_{W^n}$.
    \item For any lc center $L$ of $(W^n,B_{W^n},\Mm^{W^n})$, $\pi(L)$ is an  lc center of $(X,B,\Mm)$.
    \item For any lc center $C$ of $(X,B,\Mm)$, $\pi^{-1}(C)$ is a union of lc centers of  $(W^n,B_{W^n},\Mm^{W^n})$.
\end{enumerate}
\end{lem}
\begin{proof}
We may assume that $W$ is irreducible and $W\not=X$.

Let $f: Y\rightarrow X$ be a dlt modification (cf. \cite[Proposition 3.10]{HL22}) of $(X,B,\Mm)$, such that there exists a prime divisor $S\subset \lfloor B_Y\rfloor$ such that $f(S)=W$, where $K_Y+B_Y+\Mm_Y:=f^*(K_X+B+\Mm_X)$. Let $W_Y$ be an lc center of $(Y,B_Y,\Mm)$ which is minimal with respect to inclusion under the condition $f(W_Y)=W$. Since $(Y,B_Y,\Mm)$ is dlt, by repeatedly applying adjunction (cf. \cite[Proposition 2.10]{HL22}), we get a dlt g-pair $(W_Y,B_{W_Y},\Mm^{W_Y})/U$ such that
$$K_{W_Y}+B_{W_Y}+\Mm^{W_Y}_{W_Y}:=(K_Y+B_Y+\Mm_Y)|_{W_Y}.$$
By construction, there exists a naturally induced projective surjective morphism $f_W: W_Y\rightarrow W^n$ such that $K_{W_Y}+B_{W_Y}+\Mm^{W_Y}_{W_Y}\sim_{\Rr,W^n}0$. By \cite[Lemma 3.19]{LX22}, \cite[Theorem 2.14]{LX23}, there exists an lc g-pair $(W^n,B_{W^n},\Mm^{W^n})/U$, such that
\begin{itemize}
\item $(W^n,B_{W^n},\Mm^{W^n})$ is induced by a canonical bundle formula of $(W_Y,B_{W_Y},\Mm^{W_Y})\rightarrow W^n$,
    \item any lc center of $(W^n,B_{W^n},\Mm^{W^n})$ is the image of an lc center of  $(W_Y,B_{W_Y},\Mm^{W_Y})$, and
    \item the image of any lc center of $(W_Y,B_{W_Y},\Mm^{W_Y})$ on $W^n$ is an lc center of $(W^n,B_{W^n},\Mm^{W^n})$.
\end{itemize}
We show that  $(W^n,B_{W^n},\Mm^{W^n})$ satisfies our requirement. 

(1) It immediately follows from our construction. 

(2) $L$ is the image of an lc center $L_Y$ of $(W_Y,B_{W_Y},\Mm^{W_Y})$. By repeatedly applying \cite[Lemma 3.18]{LX22}, $L_Y$ is an lc center of $(Y,B_Y,\Mm)$. Since $K_Y+B_Y+\Mm_Y:=f^*(K_X+B+\Mm_X)$, $f(L_Y)=\pi(L)$ is an lc center of $(X,B,\Mm)$.

(3) $f^{-1}(C)$ is a union of lc centers of $(Y,B_Y,\Mm)$. Since $(Y,B_Y,\Mm)$ is dlt, $f^{-1}(C)\cap W_Y$ is a union of lc centers of $(Y,B_Y,\Mm)$. By \cite[Lemma 3.18]{LX22},  $f^{-1}(C)\cap W_Y$ is a union of lc centers of $(W_Y,B_{W_Y},\Mm^{W_Y})$. Hence $\pi^{-1}(C)=f_W(f^{-1}(C)\cap W_Y)$ is a union of lc centers of $(W^n,B_{W^n},\Mm^{W^n})$.
\end{proof}

\begin{lem}\label{lem:pushout2}
Let $(X,B,\Mm)$ be an lc g-pair. Let $Y$ and $Y'$ be two unions of lc centers, such that $Y'\subsetneq Y$, and $Y$ and $Y'$ are adjacent in $S(X,B,\Mm)$. Let $\pi: Y^n\rightarrow Y$ be the normalization of $Y$ and let $Y'':=\pi^{-1}(Y')$ with the reduced scheme structure. Denote the induced morphism $Y''\rightarrow Y'$ by $\pi''$. Then there exist a universal push-out diagram
\begin{center}
$\xymatrix{
Y'' \ar@{^(->}[r]^j\ar@{->}[d]_{\pi''} & Y^n\ar@{->}[d]^{\pi}\\
Y' \ar@{^(->}[r]^i& Y\\
}$
\end{center}
and a short exact sequence
\begin{align*}
0\to \Oo_{Y}\xrightarrow{\pi^*\oplus i^*} \pi_*\Oo_{Y^n}\oplus\Oo_{Y'}\xrightarrow{j^*-\pi''^*} \pi''_*\Oo_{Y''}\to 0,
\end{align*}
where $i,j$ are the natural closed immersions.
\end{lem}
\begin{proof}
By \cite[Theorem 4.10]{LX22} and \cite[Theorem 9.26]{Kol13}, $Y$ is semi-normal. Let $L$ be an lc center contained in $Y$ but not contained in $Y'$. Since $Y'$ and $Y$ are adjacent in $S(X,B,\Mm)$, we have
$$Y\backslash Y'=L\backslash (L\cap Y'),$$
and $L\cap Y'$ is the union of all lc centers of $(X,B,\Mm)$ that are contained in $L$ but not equal to $L$. By \cite[Theorem 4.10]{LX22}, $Y\setminus Y'$ is normal. The lemma follows from Lemma \ref{lem:pushout}.
\end{proof}

\section{Proof of the vanishing theorems}

In this section, we prove Theorem \ref{thm: kod vanishing with lc strata}, which immediately implies Theorems \ref{thm: kod vanishing gpair intro} and \ref{thm: kv vanishing gpair intro}.

\begin{lem}\label{lem:3.1}
Let $(X,B,\Mm)/U$ be an lc g-pair associated with morphism $f: X\rightarrow U$, and $D$ a Cartier divisor on $X$ such that $D-(K_X+B+\Mm_X)$ is nef$/U$ and log big$/U$ with respect to $(X,B,\Mm)$. Let $W=\Nklt(X,B,\Mm)$ with the reduced scheme structure, and let $\mathcal{I}_W$ be the defining ideal sheaf of $W$ on $X$. Then:
\begin{enumerate}
	\item $R^if_*(\mathcal{I}_W\otimes\mathcal{O}_X(D))=0$ for any $i>0$.
    \item $f_*\mathcal{O}_X(D)\rightarrow f_*\mathcal{O}_W(D)$ is surjective.
\end{enumerate}
\end{lem}
\begin{proof}
By \cite[Lemma 2.4]{Xie22}, there exists a pair $(X,\Delta)$ such that $L-K_X-\Delta$ is ample/$U$ and $W=\Nlc(X,\Delta)$. (1) follows from  \cite[Theorem 8.1]{Fuj11}. (2) follows from (1) and the long exact sequence
$$0\rightarrow f_*(\mathcal{I}_W\otimes\mathcal{O}_X(D))\rightarrow f_*\mathcal{O}_X(D)\rightarrow f_*\mathcal{O}_W(D)\rightarrow R^1f_*(\mathcal{I}_W\otimes\mathcal{O}_X(D))\rightarrow\dots.$$
\end{proof}

\begin{lem}\label{lem:3.2}
Let $(X,B,\Mm)/U$ be an lc g-pair associated with morphism $f: X\rightarrow U$, and $D$ a Cartier divisor on $X$ such that $D-(K_X+B+\Mm_X)$ is nef$/U$ and log big$/U$ with respect to $(X,B,\Mm)$. Let $Y$ and $Y'$ be two unions of lc centers, such that $Y'\subsetneq Y$, and $Y$ and $Y'$ are adjacent in $S(X,B,\Mm)$. Let $\pi: Y^n\rightarrow Y$ be the normalization of $Y$, $Y'':=\pi^{-1}(Y')$ with the reduced scheme structure, and $\pi'':=\pi|_{Y''}$.
\begin{center}
$\xymatrix{
Y'' \ar@{^(->}[r]^j\ar@{->}[d]_{\pi''} & Y^n\ar@{->}[d]^{\pi}\\
Y' \ar@{^(->}[r]^i& Y\\
}$
\end{center}
Then the induced map
$$f_*\pi_*\mathcal{O}_{Y^n}(D|_{Y^n}) \rightarrow f_*\pi''_*\mathcal{O}_{Y''}(D|_{Y''})$$
is surjective.
\end{lem}

\begin{proof}
By Lemma \ref{lem:adjunction}, there exists an lc g-pair $(Y^n,B_{Y^n},\Mm^{Y^n})/U$, such that 
\begin{itemize}
    \item $K_{Y^n}+B_{Y^n}+\Mm^{Y^n}_{Y^n}\sim_{\Rr,U}(K_X+B+\Mm_X)|_{Y^n}$,
    \item for any lc center $L$ of $(Y^n,B_{Y^n},\Mm^{Y^n})$, $\pi(L)$ is an  lc center of $(X,B,\Mm)$, and
    \item for any lc center $C$ of $(X,B,\Mm)$, $\pi^{-1}(C\cap Y)$ is a union of lc centers of  $(Y^n,B_{Y^n},\Mm^{Y^n})$.
\end{itemize}
Then $D|_{Y^n}-(K_{Y^n}+B_{Y^n}+\Mm^{Y^n}_{Y^n})$ is nef$/U$ and log big$/U$ with respect to $(Y^n,B_{Y^n},\Mm^{Y^n})$.

Let $Y_0$ be a connected component of $Y^n$ and let $Y_0'':=Y''\cap Y_0$. We claim that 
\begin{align}\label{eq}
\text{either $Y_0''=Y_0$ or $Y_0''=\Nklt(Y_0,B_{Y^n}|_{Y^0},\Mm^{Y^n}|_{Y^0})$.} 
\end{align}
Indeed, if this is not the case, then there exists an lc center $L$ of $(Y_0,B_{Y^n}|_{Y^0},\Mm^{Y^n}|_{Y^0})$ such that $L\neq Y_0$ and $L$ is not contained in $Y_0''$.
Then $\tilde{Y}:=\pi(Y''\cup L)\in S(X,\Delta,\Mm)$ and $Y' \subsetneq \tilde{Y}\subsetneq Y$, which contradicts the condition that $Y',Y$ are adjacent. By \eqref{eq} and Lemma \ref{lem:3.1}(2), 
$$f_*\pi_*\mathcal{O}_{Y_0}(D|_{Y_0}) \rightarrow f_*\pi''_*\mathcal{O}_{Y_0''}(D|_{Y_0''})$$
is surjective. Thus 
$$f_*\pi_*\mathcal{O}_{Y^n}(D|_{Y^n}) \rightarrow f_*\pi''_*\mathcal{O}_{Y''}(D|_{Y''})$$
is surjective.
\end{proof}

\begin{proof}[Proof of Theorem \ref{thm: kod vanishing with lc strata}]
We apply induction on $\dim X$. When $\dim X=1$ the theorem is obvious. 

For any union of lc centers $Z$ of $(X,B,\Mm)$, we define $m(Z)$ to be the number of lc centers of $(X,B,\Mm)$ that are contained in $Z$. We let $W:=\Nklt(X,B,\Mm)$, associated with the reduced scheme structure.

\medskip

\noindent\textbf{Step 1}. In this step we prove (1) when $Y$ is minimal in $S(X,B,\Mm)$.

By \cite[Theorem 4.10]{LX22}, $Y$ is normal. If $\dim Y=0$ then we are done. Otherwise, by Lemma \ref{lem:adjunction}, there exists a klt g-pair $(Y,B_Y,\Mm^{Y})/U$ such that $K_{Y}+B_{Y}+\Mm^{Y}_{Y}\sim_{\Rr,U}(K_X+B+\Mm_X)|_{Y}$. Hence $D|_Y-(K_{Y}+B_{Y}+\Mm^{Y}_{Y})$ is nef$/U$ and big$/U$. By \cite[Lemma 2.4]{Xie22}, there exists a klt pair $(Y,\Delta_Y)$ such that $D|_Y-(K_{Y}+\Delta_Y)$ is ample$/U$. (1) follows from the usual Kawamata-Viehweg vanishing theorem (cf. \cite[Theorem 1-2-7]{KMM87}).

\medskip

\noindent\textbf{Step 2}. In this step we prove (1). 

We apply induction on $m(Y)$. When $m(Y)=1$, $Y$ is minimal in $S(X,B,\Mm)$ and we are done by \textbf{Step 1}. Thus we may assume that $m(Y)>1$. Then there exists a union of lc centers $Y'$ such that $Y'\subsetneq Y$, and $Y$ and $Y'$ are adjacent in $S(X,B,\Mm)$. Since $m(Y')<m(Y)$, by induction on $m(Y)$, we have
\begin{align}\label{eq3.1}
R^if_*\mathcal{O}_{Y'}(D)=0
\end{align}
for any positive integer $i$. 

Let $\pi: Y^n\rightarrow Y$ be the normalization of $Y$, and let $Y'':=\pi^{-1}(Y')$ with the reduced scheme structure. Let $i: Y'\hookrightarrow Y$ and  $j: Y''\hookrightarrow Y^n$  be the natural inclusions, and let $\pi'':=\pi|_{Y''}$. By Lemma \ref{lem:pushout2}, there exists a universal push-out diagram
\begin{center}
$\xymatrix{
Y'' \ar@{^(->}[r]^j\ar@{->}[d]_{\pi''} & Y^n\ar@{->}[d]^{\pi}\\
Y' \ar@{^(->}[r]^i& Y\\
}$
\end{center}
and a short exact sequence
\begin{align}\label{eq: short exact sequence in main theorem}
0\to \Oo_{Y}\xrightarrow{\pi^*\oplus i^*} \pi_*\Oo_{Y^n}\oplus\Oo_{Y'}\xrightarrow{j^*-\pi''^*} \pi''_*\Oo_{Y''}\to 0.
\end{align}
By Lemma \ref{lem:adjunction}, there exists an lc g-pair $(Y^n,B_{Y^n},\Mm^{Y^n})/U$, such that 
\begin{itemize}
    \item $K_{Y^n}+B_{Y^n}+\Mm^{Y^n}_{Y^n}\sim_{\Rr,U}(K_X+B+\Mm_X)|_{Y^n}$,
    \item for any lc center $L$ of $(Y^n,B_{Y^n},\Mm^{Y^n})$, $\pi(L)$ is an  lc center of $(X,B,\Mm)$, and
    \item for any lc center $C$ of $(X,B,\Mm)$, $\pi^{-1}(C\cap Y)$ is a union of lc centers of  $(Y^n,B_{Y^n},\Mm^{Y^n})$.
\end{itemize}
Then $D|_{Y^n}-(K_{Y^n}+\Delta_{Y^n}+\Mm^{Y^n}_{Y^n})$ is nef$/U$ and log big$/U$ with respect to $(Y^n,B_{Y^n},\Mm^{Y^n})$, and $Y''$ is a union of lc centers of $(Y^n,B_{Y^n},\Mm^{Y^n})$. Since $\dim Y^n<\dim X$ and $\pi$ is a finite morphism, by  induction on $\dim X$ we have
\begin{align}\label{eq3.2}
&R^i(f\circ \pi)_*\mathcal{O}_{Y^n}(D|_{Y^n})=R^if_*\big(\pi_*(\mathcal{O}_{Y^n}(D|_{Y^n})\big)
=0
\end{align}
and 
\begin{align}\label{eq3.3}
&R^i(f\circ \pi'')_*\mathcal{O}_{Y''}(D|_{Y''})= R^if_*\big(\pi''_*\mathcal{O}_{Y''}(D|_{Y''})\big)=0.
\end{align}
By the short exact sequence \eqref{eq: short exact sequence in main theorem}, we have a short exact sequence
$$0\rightarrow \mathcal{O}_Y(D)\xrightarrow{\pi^*\oplus i^*}\pi_*\mathcal{O}_{Y^n}(D|_{Y^n})\oplus \mathcal{O}_{Y'}(D) \xrightarrow{j^*-\pi''^*} \pi''_*\mathcal{O}_{Y''}(D|_{Y''})\rightarrow 0,$$
which induces a long exact sequence
\begin{align*}
0&\rightarrow f_*\mathcal{O}_Y(D)\rightarrow f_*\pi_*\mathcal{O}_{Y^n}(D|_{Y^n})\oplus f_*\mathcal{O}_{Y'}(D)  \xrightarrow{j^*-\pi''^*} f_*\pi''_*\mathcal{O}_{Y''}(D|_{Y''})\rightarrow \cdots\\
\cdots &\rightarrow R^if_*\mathcal{O}_Y(D) \rightarrow R^if_*\big(\pi_*(\mathcal{O}_{Y^n}(D|_{Y^n})\big)\oplus R^if_*\mathcal{O}_{Y'}(D)\rightarrow  R^if_*\big(\pi''_*\mathcal{O}_{Y''}(D|_{Y''})\big)\rightarrow\cdots.
\end{align*}
Hence, it follows from \eqref{eq3.1}, \eqref{eq3.2}, \eqref{eq3.3} and Lemma \ref{lem:3.2} that $R^if_*\mathcal{O}_Y(D)=0$ for any positive integer $i$.

\medskip

\noindent\textbf{Step 3}. In this step we prove (2) and prove (3)(4) when $Y=W=\Nklt(X,B,\Mm)$. 

We have the long exact sequence
\begin{align*}
    0&\rightarrow f_*(\mathcal{I}_W\otimes \mathcal{O}_X(D))\rightarrow f_*\mathcal{O}_X(D)\rightarrow f_*\mathcal{O}_W(D)\rightarrow\dots\\
    \dots&\rightarrow R^if_*(\mathcal{I}_W\otimes \mathcal{O}_X(D))\rightarrow R^if_*\mathcal{O}_X(D)\rightarrow R^if_*\mathcal{O}_W(D)\rightarrow\dots
\end{align*}
By (1), $R^if_*\mathcal{O}_W(D)=0$ for any positive integer $i$. By Lemma \ref{lem:3.1}(1), $R^i(\mathcal{I}_W\otimes f_*\mathcal{O}_X(D))=0$ for any positive integer $i$. This implies (2), and also implies (3)(4) when $Y=W$.

\medskip

\noindent\textbf{Step 4}. We prove (3)(4) in this step, hence conclude the proof of the theorem.

We apply induction on $m(W)-m(Y)$. When $m(W)-m(Y)=0$, $Y=W$ and we are done by \textbf{Step 3}. Thus we may assume that $m(W)-m(Y)>0$. Then there exists a union of lc centers $\tilde Y$ such that $Y\subsetneq\tilde Y\subset W$, and $Y$ and $\tilde Y$ are adjacent in $S(X,B,\Mm)$.

Let $\tilde\pi:\tilde{Y}^n\rightarrow \tilde{Y}$ be the normalization of $\tilde Y$, and let $\widehat Y:=\tilde\pi^{-1}(Y)$ with the reduced scheme structure. Let $\tilde i: Y\hookrightarrow\tilde Y$ and  $\tilde j: \widehat{Y}\hookrightarrow\tilde Y^n$  be the natural inclusions, and let $\widehat{\pi}:=\tilde\pi|_{\widehat{Y}}$. By Lemma \ref{lem:pushout2}, there exists a universal push-out diagram
\begin{center}
$\xymatrix{
\widehat{Y}\ar@{^(->}[r]^{\tilde j}\ar@{->}[d]_{\widehat\pi} & \tilde Y^n\ar@{->}[d]^{\tilde \pi}\\
Y\ar@{^(->}[r]^{\tilde i}&\tilde Y\\
}$
\end{center}
and a short exact sequence
\begin{align*}
0\to \Oo_{\tilde Y}\xrightarrow{\tilde\pi^*\oplus\tilde i^*} \tilde\pi_*\Oo_{\tilde Y^n}\oplus\Oo_{Y}\xrightarrow{\tilde j^*-\widehat{\pi}^*} \widehat{\pi}_*\Oo_{\widehat{Y}}\to 0.
\end{align*}
which induces a short exact sequence
$$0\rightarrow \mathcal{O}_{\tilde{Y}}(D)\xrightarrow{\tilde\pi^*\oplus\tilde i^*} \tilde\pi_*\mathcal{O}_{\tilde{Y}^n}(D|_{\tilde{Y}^n})\oplus \mathcal{O}_{Y}(D)  \xrightarrow{\tilde j^*-\widehat{\pi}^*} \widehat{\pi}_*\mathcal{O}_{\widehat{Y}}(D|_{\widehat{Y}})\rightarrow 0.$$
So we have the left exact sequence
\begin{align}\label{eq:long}
0\rightarrow f_*\mathcal{O}_{\tilde{Y}}(D)\xrightarrow{\tilde\pi^*\oplus\tilde i^*} f_*\tilde\pi_*\mathcal{O}_{\tilde{Y}^n}(D|_{\tilde{Y}^n})\oplus f_*\mathcal{O}_{Y}(D)  \xrightarrow{\tilde j^*-\widehat{\pi}^*} f_*\widehat{\pi}_*\mathcal{O}_{\widehat{Y}}(D|_{\widehat{Y}}).
\end{align}
By Lemma \ref{lem:3.2},
$$\tilde j^*: f_*\tilde\pi_*\mathcal{O}_{\tilde{Y}^n}(D|_{\tilde{Y}^n})\rightarrow f_*\widehat{\pi}_*\mathcal{O}_{\widehat{Y}}(D|_{\widehat{Y}})$$
is surjective. Thus by an easy map tracing of \eqref{eq:long} we have that
$$\tilde i^*: f_*\mathcal{O}_{\tilde{Y}}(D)\rightarrow f_*\mathcal{O}_{Y}(D)$$
is also surjective. Since  $m(W)-m(\tilde Y)<m(W)-m(Y)$, by induction on $m(W)-m(Y)$, $$f_*\mathcal{O}_X(D)\rightarrow f_*\mathcal{O}_{ \tilde{Y}}(D)$$ is surjective. This implies (3).

We have the long exact sequence
\begin{align*}
    0&\rightarrow f_*(\mathcal{I}_Y\otimes \mathcal{O}_X(D))\rightarrow f_*\mathcal{O}_X(D)\rightarrow f_*\mathcal{O}_Y(D)\rightarrow\dots\\
    \dots&\rightarrow R^if_*(\mathcal{I}_Y\otimes \mathcal{O}_X(D))\rightarrow R^if_*\mathcal{O}_X(D)\rightarrow R^if_*\mathcal{O}_Y(D)\rightarrow\dots,
\end{align*}
so (4) follows immediately from (1)(2)(3).
\end{proof}

\begin{proof}[Proof of Theorem \ref{thm: kod vanishing gpair intro}]
It immediately follows from Theorem \ref{thm: kod vanishing with lc strata}(2) by letting $U=\{pt\}$.
\end{proof}

\begin{proof}[Proof of Theorem \ref{thm: kv vanishing gpair intro}]
It immediately follows from Theorem \ref{thm: kod vanishing with lc strata}(2).
\end{proof}

\section{Base-point-freeness for lc g-pairs}

In this section, we prove Theorems \ref{thm:base-point-freeness intro}, \ref{thm: semi-ampleness intro}, and \ref{thm: cont thm gpair}.

\begin{lem}\label{lem: non-vanishing of lc gpair}
 Let $a$ be a positive real number, $(X,B,\Mm)/U$ an lc g-pair, and $D$ a nef$/U$ Cartier divisor on $X$ such that $aD-(K_X+B+\Mm_X)$ is ample$/U$. Let $Y$ be a minimal lc center of $(X,B,\Mm)$ if $(X,B,\Mm)$ is not klt, and let $Y:=X$ if $(X,B,\Mm)$ is klt. Let $D_Y:=D|_Y$. Then for any integer $m\gg 0$,
\begin{enumerate}
    \item  $\mathcal{O}_{Y}(mD_Y)$ is globally generated over $U$,
    \item  $|mD/U|\not=\emptyset$, and
    \item $Y$ is not contained in $\Bs|mD/U|$.
\end{enumerate}
\end{lem}
\begin{proof}
When $(X,B,\Mm)$ is klt, by \cite[Lemma 2.4]{Xie22}, there exists a klt pair $(X,\Delta)$ such that $D-(K_X+\Delta)$ is ample$/U$. By the usual base-point-freeness theorem (cf. \cite[Theorem 3-1-1]{KMM87}), the lemma follows.

When $(X,B,\Mm)$ is not klt, by \cite[Theorem 4.10]{LX22}, $Y$ is normal. By Theorem \ref{thm: kod vanishing with lc strata}(3), the map $f_*\mathcal{O}_X(mD)\rightarrow f_*\mathcal{O}_Y(mD_Y)$ is surjective for any positive integer $m\geq a$. Thus (2)(3) follow from (1) and we only need to prove (1). If $\dim Y=0$ then there is nothing left to prove. If $\dim Y>0$, then by Lemma \ref{lem:adjunction}, there exists a klt g-pair $(Y,B_Y,\Mm^{Y})/U$ such that $K_{Y}+B_{Y}+\Mm^{Y}_{Y}\sim_{\Rr,U}(K_X+B+\Mm_X)|_{Y}$ and $\Nklt(Y,B_Y,\Mm^Y)=\Nklt(X,B,\Mm)|_Y$. Thus $D_Y-(K_{Y}+B_{Y}+\Mm^{Y}_{Y})$ is nef$/U$ and big$/U$ with respect to $(Y,B_Y,\Mm^Y)$. By \cite[Lemma 2.4]{Xie22}, there exists a klt pair $(Y,\Delta_Y)$ such that $D_Y-(K_Y+\Delta_Y)$ is ample$/U$. By the usual base-point-freeness theorem (cf. \cite[Theorem 3-1-1]{KMM87}), the lemma follows.
\end{proof}

\begin{proof}[Proof of Theorem \ref{thm:base-point-freeness intro}]
By Lemma \ref{lem: non-vanishing of lc gpair}, we may let $m_0$ be the minimal positive integer such that $|mD|\not=\emptyset$ for any integer $m\geq m_0$.

\begin{claim}\label{claim: induction bs}
Let $\{p_i\}_{i=1}^{+\infty}$ be a strictly increasing sequence of positive integers. There exist a non-negative integer $M$ and integers $i_1<i_2<\dots<i_{M+1}$ satisfying the following. Let $s_k:=\prod_{l=1}^kp_{i_l}$ for any $1\leq k\leq M+1$, then
\begin{enumerate}
    \item $|s_1D/U|\not=\emptyset$,
    \item $\Bs|s_kD/U|\supsetneq\Bs|s_{k+1}D/U|$ for any $1\leq k\leq M$, and
    \item $\Bs|s_{M+1}D/U|=\emptyset$.
\end{enumerate}
\end{claim}
\begin{proof}
We may take $i_1$ to be any integer such that $p_{i_1}\geq m_0$, then (1) holds. 

Suppose that we have already found $i_1,\dots,i_k$ for some positive integer $k$. Let $d:=\dim X$, let $H_{1},\cdots,H_{d+1}$ be $d+1$ be general elements in $|s_kD/U|$, and let $H:=H_{1}+\cdots+H_{d+1}$. Then $(X,B+H,\Mm)$ is lc outside $\Bs|s_kD/U|$. If $\Bs|s_kD/U|=\emptyset$, then we may let $M:=k-1$ and we are done. Thus we may assume that $\Bs|s_kD/U|\not=\emptyset$.

Since every $H_{j}$ contains $\Bs|s_kD/U|$, by \cite[Theorem 18.22]{Kol+92}, $(X,B+H,\Mm)$ is not lc near $\Bs|s_kD/U|$. Let 
$$c:=\sup\{t\mid t\geq 0, (X,B+tH,\Mm)\text{ is lc}\},$$
then $c\in [0,1)$, and there exists at least one lc center of $(X,B+cH,\Mm)$ which is contained in $\Bs|s_kD/U|$. Let
$\mathcal{S}$ be the set of all lc centers of $(X,B+cH,\Mm)$ that are contained in $\Bs|s_kD/U|$, and let $Y$ be a minimal lc center in $\mathcal{S}$. Since
$$(a+s_k(d+1))D-(K_X+B+cH+\Mm_X)\sim_{\mathbb R}s_k(d+1)(1-c)D+(aD-(K_X+B+\Mm_X))$$
is ample$/U$, by Lemma \ref{lem: non-vanishing of lc gpair}, there exists a positive integer $n$, such that for any integer $m\geq n$, $|ms_kD/U|\not=\emptyset$ and $\Bs|ms_kD/U|$ does not contain $Y$. In particular, $\Bs|ms_kD/U|\subsetneq\Bs|s_kD/U|$. We may let $i_{k+1}$ be any integer such that $i_{k+1}>i_{k}$ and $p_{i_{k+1}}\geq n$. This construction implies (2). (3) follows from (2) and the Noetherian property.
\end{proof}

\noindent\textit{Proof of Theorem \ref{thm:base-point-freeness intro} continued}. We let $\{p_i\}_{i=1}^{+\infty}$ and $\{q_j\}_{j=1}^{+\infty}$ be two strictly increasing sequence of prime numbers, such that $p_i\not=q_j$ for any $i,j$. By Claim \ref{claim: induction bs}, there exist two non-negative integers $M,N$ and positive integers $i_1<i_2<\dots<i_{M+1}$ and $j_1<j_2<\dots<j_{N+1}$, such that $\mathcal{O}_X(\prod_{l=1}^{M+1}p_{i_l}D)$ and $\mathcal{O}_X(\prod_{l=1}^{N+1}q_{i_l}D)$ are globally generated$/U$. Let $p:=\prod_{l=1}^{M+1}p_{i_l}$ and $q:=\prod_{l=1}^{N+1}q_{i_l}$, then $p$ and $q$ are coprime. Therefore, for any integer $m\gg 0$, we may write $m=bp+cq$ for some non-negative integers $b,c$, hence
$$\Bs|mD/U|\subset\Bs|pD/U|\cup\Bs|qD/U|=\emptyset.$$
Therefore, $\mathcal{O}_X(mD)$ is globally generated over $U$ for any integer $m\gg 0$.
\end{proof}

\begin{proof}[Proof of Theorem \ref{thm: semi-ampleness intro}]
By the theory of Shokurov-type rational polytopes (cf. \cite[Proposition 3.20]{HL22}) and the theory of uniform rational polytopes (cf. \cite[Lemma 5.3]{HLS19}, \cite[Theorem 1.4]{Che20}), we may assume that $D$ is a $\Qq$-divisor. The theorem immediately follows from Theorem \ref{thm:base-point-freeness intro}. 
\end{proof}

\begin{proof}[Proof of Theorem \ref{thm: cont thm gpair}]
(1)(2) By the cone theorem \cite[Theorem 1.1(1-4)]{HL21},  $F$ is a finitely dimensional rational  $(K_X+B+\Mm_X)$-negative extremal face$/U$. Thus there exists a nef Cartier divisor $L$ on $X$ that is the supporting function of $F$. Then $L-(K_X+B+\Mm_X)$ is ample. By Theorem \ref{thm:base-point-freeness intro}, $mL$ is base-point-free$/U$, hence defines a contraction$/U$. Denote this contraction by $\cont_F$. Then $\cont_F$ satisfies (1) and (2).

(3) Since $D-(K_X+B+\Mm_X)$ is ample$/Z$, by Theorem \ref{thm:base-point-freeness intro}, $\mathcal{O}_X(mD)$ is globally generated over $Z$ for any integer $m\gg 0$. Since $D\cdot C$ for any curve $C$ contracted by $\cont_F$, $\cont_F$ is defined by $|mD|$ for any integer $m\gg 0$. Thus $mD=f^*D_{Y,m}$ and $(m+1)D=f^*D_{Y,m+1}$ for any integer $m\gg 0$. We may let $D_Y:=D_{Y,m+1}-D_{Y,m}$.
\end{proof}

\begin{rem}\label{rem: does not rely on kollar}
Koll\'ar's gluing theory for generalized pairs was originally established in \cite[Construction 4.12]{LX22} to glue glc crepant structures. This theory was further developed in \cite{Xie22}. Although we have extensively referenced both \cite{LX22} and \cite{Xie22}, it is important to note that we have only cited results before \cite[Theorem 4.10]{LX22} from \cite{LX22} and only cited \cite[Lemma 2.4]{Xie22} from \cite{Xie22}. None of these cited results are dependent on Koll\'ar's gluing theory for generalized pairs (although \cite[Theorem 4.10]{LX22} used the idea of stratification). Consequently, the proofs of our main theorems are independent of Koll\'ar's gluing theory for generalized pairs.
\end{rem}

\end{document}